\documentclass{article}
\usepackage{amsmath}
\usepackage{amsfonts}
\usepackage{amssymb}
\usepackage{euscript}

\newtheorem{thm}{Theorem}[section]
\newtheorem{cor}{Corollary}[section]

\newtheorem{example}{Example}[section]
\newtheorem{remark}{Remark}[section]

\def\dj{d\kern-0.4em\char"16\kern-0.1em}

\title{
Finite and infinite order differential properties\\
of the reduced Mittag--Leffler polynomials }
\author{
{\sc Predrag M. Rajkovi\'c, }\\
[2mm]
{\sc Sladjana D. Marinkovi\'c,}\\
[2mm]
{\sc Miomir S. Stankovi\'c,}\\
[2mm]
{\sc Marko D. Petkovi\'c}\\
[2mm] {\bf   University of Ni\v s,\ Serbia} }

\date{}

\begin{document}
\maketitle

\medskip

{\it Abstract.} This paper deals with the Mittag-Leffler
polynomials (MLP) by extracting their essence which consists of
real polynomials with fine properties. They are orthogonal on the
real line instead of the imaginary axes for MLP. Beside recurrence
relations and zeros, we will point to the closed form of its
Fourier transform. The most important contribution consists of the
new differential properties, especially the finite and infinite
differential equation.

\medskip

% 33C45 Orthogonal polynomials and functions of hypergeometric type
%(Jacobi, Laguerre, Hermite, % Askey scheme, etc.)
% 11B83 Special sequences and polynomials
% 68W30 Symbolic computation and algebraic computation
% 40A30 Convergence and divergence of series and
%sequences of functions

% 30C15 Zeros of polynomials, rational functions, and other analytic
%     functions
%42C05 Orthogonal functions and polynomials, general theory
%%% 33E12 Mittag-Leffler functions and generalizations



\section{Introduction}

The Mittag--Leffler polynomials $\{g_n(x)\}$ are coefficients in
the expansion
$$
\Bigl(\frac{1+t}{1-t}\Bigr)^{x}=\sum_{n=0}^\infty g_n(x)t^n\qquad
(|t|<1).
$$
They were introduced by Mittag-Leffler in a study of the integral
representations. Their main properties were discovered by H.
Bateman (see \cite{Bateman1} and \cite{Bateman2}). He noticed that
they occur as coefficients in the closed-form expressions for a
several families of integrals. Also, they were used in deriving
some expansions for the Euler Gamma function and the Riemann Zeta
function \cite{Rzadkowski}. Truncated Exponential-Based
Mittag-Leffler Polynomials were examined in \cite{Yasmin}. They
were noticed in the solutions of heat diffusion equation of
Fokker-Plank type and researched by G. Dattoli and his coauthors
\cite{Gorska1}-\cite{Gorska2}. They are connected with the Sheffer
polynomials in \cite{He} and Riordan arrays in \cite{Luzon}. Their
generalizations were considered in \cite{StaMarRaj} and
\cite{Raza} .

The article is organized as follows. In Section 2 we present the
preliminaries for the Mittag-Leffler polynomials. Since the
imaginary unit is present in the orthogonality relation, we have
noticed that we can reduce them to the real polynomials whose
examination is more obvious and much easier. They are the subject
in Section 3. In the next section, we find the closed form of
their Fourier transform. The main contributions of the paper are
the differential properties exposed Section 5. Finally, the paper
is concluded in Section 6 with observation on the quasi-monomial
property of the reduced Mittag-Leffler polynomials.

\section{Preliminaries}

The Mittag--Leffler polynomials $\{g_n(x)\}$  can be represented
over hypergeometric function like
\begin{equation}\label{hyper}
g_n(x)=2x\;{}_2F_1\Bigl({{1-n,\ 1-x}\atop{2}}\Bigm| 2\Bigr)\qquad
(n\in\mathbb N).
\end{equation}
They can be viewed a special case of the Meixner-Pollaczek
polynomials
\begin{equation}\label{Meixner-Pollaczek-1}
P^{(\lambda)}_n(x;\phi)=\frac{(2\lambda)_n}{n!}\,{\rm e}^{{\rm
i}n\phi} \;{}_2F_1\Bigl({{-n,\ \lambda+{\rm
i}x}\atop{2\lambda}}\Bigm| 1-{\rm e}^{-2{\rm i}\phi}\Bigr)\quad
(n\in\mathbb N;\ {\rm i}^2=-1)
\end{equation}
for $\lambda = 1$ and $\phi = \pi/2$:% and $y={\rm i}x$:
\begin{equation}\label{Meixner-Pollaczek-2}
g_n(x)=2\frac{{\rm e}^{-{\rm i}n\pi/2}}{n}\,xP^{(1)}_{n-1}({\rm
i}x;\pi/2)\qquad (n\in\mathbb N).
\end{equation}
Also, they can be considered as a special case of Meixner
polynomials
%\begin{equation}\label{Meixner}
$$
M_n(x;\beta,c) = \;{}_2F_1\Bigl({{-n,\ -x}\atop{\beta}}\Bigm|
1-\frac1{c}\Bigr)\qquad (n\in\mathbb N)
$$
%\end{equation}
for $\beta=2$ and $c=-1$ (see \cite{Koek}). Namely, their relation
is given by
$$
g_n(x)=2xM_{n-1}(x-1,2,-1).
$$
The lack of this connection is the fact that orthogonality of the
Meixner polynomials is assured only with the constraint $0<c<1$.

Finally, the Mittag--Leffler polynomials $g_n(x)$ are connected
with Pidduck polynomials \cite{Weisstein} by the expression
$$
P_n(x)=\dfrac{1}{2}\left({\rm e}^D+1\right)g_n(x),
$$
where the series for the exponential function is used and $D$ is
understood as differentiation.

The Mittag--Leffler polynomials $\{g_n(x)\}$ satisfy recurrence
relation
\begin{equation}\label{MLPRec-1}
(n+1)g_{n+1}(x)-2xg_n(x)+(n-1)g_{n-1}(x)=0,
\end{equation}
with initial values
\begin{equation*}%\label{MLPRec}
g_0(x)=1,\quad g_1(x)=2x.
\end{equation*}
They also % Mittag--Leffler polynomials $\{g_n(x)\}$
satisfy difference relation
\begin{equation*}%\label{MLPRec-2}
xg_{n}(x+1)-2ng_n(x)-xg_{n}(x-1)=0.
\end{equation*}
The Mittag--Leffler polynomials $\{g_n(x)\}$  satisfy
recurrence-difference  relation
\begin{equation*}%\label{MLPRec-3}
g_n(x+1)-g_{n-1}(x+1)=g_n(x)+g_{n-1}(x).
\end{equation*}
The orthogonality relation is given by
\begin{equation}\label{MLPOrtho}
\int_{-\infty}^{+\infty}g_n(-{\rm i}x)g_m({\rm i}x)\frac{{\rm
d}x}{x\sinh (\pi x)}=\frac{2}{n}\delta_{mn}\qquad(n,m\in\mathbb
N).
\end{equation}
Notice that the corresponding monic sequence are
\begin{equation*}%\label{MLPMonic}
\hat g_n(x) = \frac{(n)!}{2^{n}}\ g_n(x)  \qquad (n\in\mathbb
N_0).
\end{equation*}
Let remind that the central difference operator is
\begin{equation*}%\label{Delta}
\delta f(z) = f(z+1/2) - f(z-1/2).
\end{equation*}
\begin{thm}
The Mittag--Leffler polynomials $\{g_n(x)\}$  satisfy the
Rodrigues formula
\begin{equation*}%\label{Rodr}
g_n(x) = \frac{2}{n!} \frac{x}{w(x,1)}\delta^n w(x,n),
\end{equation*}
where
\begin{equation*}%\label{Weight}
w(x,n) = \Gamma\left(\frac{n+1}2-x\right)
\Gamma\left(\frac{n+1}2+x\right).
\end{equation*}
\end{thm}
{\it Proof}. It is based on the connection (\ref{Meixner-Pollaczek-2})
and the Rodrigues formula for the Meixner-Pollaczek polynomials in \cite{Koek}, pp. 37.$\Box$

\section{The reduced Mittag--Leffler polynomials}

In spite of the fact that the Mittag-Leffler polynomials are real,
in some relations, as in (\ref{Meixner-Pollaczek-1}), and
especially in their orthogonality relation (\ref{MLPOrtho}) they
are considered as the complex functions. That is why we believe
that is better to extract the following sequence.

 Let us consider {\it reduced Mittag-Leffler polynomials}
defined by (see \cite{StaMarRaj})
\begin{equation*}%\label{ModMLP}
\varphi_n(x) = \frac{g_{n+1}({\rm i}x)}{{\rm i}^{n+1}\ x}\qquad
(n\in\mathbb N_0).
\end{equation*}

The successive members of sequence $\{\varphi_n(x)\}_{n\in\mathbb
N_0}$ satisfy the three--term recurrence relation
\begin{eqnarray*}%\label{ModMLPRec}
(n+2)\varphi_{n+1}(x) &=& 2x\varphi_{n}(x) - n\ \varphi_{n-1}(x) \quad (n\in \mathbb N)\\
\varphi_0(x)&=&2, \qquad \varphi_1(x)=2x.\nonumber
\end{eqnarray*}
The generating function of sequence $\{\varphi_n(x)\}_{n\in\mathbb
N_0}$ is given by
\begin{equation*}%\label{GenFunvarphi}
\mathcal G(t,x)=\frac{\exp(2x\arctan t)-1}{tx}=\sum_{n=0}^\infty
\varphi_n(x)t^n.
\end{equation*}
Notice that
\begin{equation*}%\label{varphi(-x)}
\varphi_n(-x)=(-1)^n \ \varphi_n(x) \qquad (n\in\mathbb N).
\end{equation*}
The polynomials of sequence $\{\varphi_n(x)\}_{n\in\mathbb N_0}$
satisfy following orthogonality relation:
\begin{equation*}%\label{ReducedMLP}
\int_{-\infty}^{+\infty}\varphi_n(x)\varphi_m(x)\frac{x}{\sinh
(\pi x)}\ {\rm d}x = \frac{2}{n+1}\
\delta_{mn}\qquad(n,m\in\mathbb N_0).
\end{equation*}
Notice that
\begin{equation*}%\label{ReducedMLP-2}
\int_{-\infty}^{+\infty}\frac{x^n\;{\rm d}x}{\sinh x} =
\frac{(1-(-1)^n)(2^{n+1}-1)}{2^n}\;n!\;\zeta(n+1)\qquad(n\in\mathbb
N),
\end{equation*}
where $\zeta(n)$ is the Riemann zeta function.

 The reduced Mittag-Leffler polynomials
$\{\varphi_n(x)\}_{n\in\mathbb N_0}$ are real polynomials and,
because of orthogonality, they have all real zeros.

The monic sequence
\begin{equation}\label{ModMLPMonic}
\hat\varphi_n(x) = \frac{(n+1)!}{2^{n+1}}\ \varphi_n(x)  \qquad
(n\in\mathbb N_0)
\end{equation}
satisfies three term recurrence relation
\begin{eqnarray}\label{ModMLPRecMonic}
\hat\varphi_{n+1}(x) &=& x\hat\varphi_{n}(x) - \frac{n(n+1)}{4}\ \hat\varphi_{n-1}(x) \quad (n\in \mathbb N),\\
\hat\varphi_0(x)&=&1, \qquad \hat\varphi_1(x)=x.\nonumber
\end{eqnarray}
\begin{remark}\rm  Notice that $\varphi_n(x)$ is the
Meixner-Pollaczeck polynomial for $\lambda=1$ and $\phi=\pi/2$,
i.e.
$$
\hat\varphi_n(x) = \frac{n!}{2^{n}}\; P_n^{(1)}(x;\pi/2) \qquad
(n\in\mathbb N_0).
$$
Hence the difference relation
$$
(x+{\rm i})\hat\varphi_n(x+{\rm i}) - 2(n+1){\rm
i}\hat\varphi_n(x) -(x-{\rm i})\hat\varphi_n(x-{\rm i})=0
$$
is valid.
\end{remark}

\begin{example}\rm
The first members of the sequence
$\{\hat\varphi_n(x)\}_{n\in\mathbb N_0}$ are
$$
\hat\varphi_{0}(x) = 1,\quad \hat\varphi_{1}(x) = x,\quad
\hat\varphi_{2}(x) = x^2-\frac12,\quad \hat\varphi_{3}(x) =
x^3-2x,
$$
$$
\hat\varphi_{4}(x) = x^4-5x^2+\frac32,\quad \hat\varphi_{5}(x) =
x^5-10x^3+\frac{23}2x.
$$
The largest zeros are
$$
x^{(2)}_2\approx 0.707,\quad x^{(3)}_3\approx 1.414,\quad
x^{(4)}_4\approx 2.163,\quad x^{(5)}_5\approx 2.945.
$$
\end{example}
Using conclusions from the paper \cite{Ismail}, we conclude that the zeros $\{x_k^{(n)}\}$ of the
polynomial $\hat\varphi_n(x)$ are bordered in the next manner:
\begin{equation*}%\label{Zeros-1}
|x_k^{(n)}| < \sqrt{(n-1)n}\qquad (k=1,2,\ldots,n).
\end{equation*}
%ILIA KRASIKOV, TURAN INEQUALITIES AND ZEROS OF ORTHOGONAL
%POLYNOMIALS, METHODS AND APPLICATIONS OF ANALYSIS, Vol. 12, No. 1,
%pp. 075--088, March 2005. Theorem 4
\begin{thm}
The sequence $\{\hat\varphi_n(x)\}_{n\in\mathbb N_0}$ satisfies
the Turan's inequality
\begin{equation*}%\label{TuranIneq}
\mathcal{T}\left(\hat\varphi_{n}(x),x\right) = - \begin{vmatrix}
\hat\varphi_{n-1}(x) & \hat\varphi_{n}(x)\\
\hat\varphi_{n}(x) & \hat\varphi_{n+1}(x)
\end{vmatrix}
\ge 0\qquad (\forall x\in\mathbb R;\ \forall n\in\mathbb N).
\end{equation*}
\end{thm}
{\it Proof}. We will prove by the mathematical induction as in the
paper \cite{Krasikov}. Obviously,
$\mathcal{T}\left(\hat\varphi_{0}(x),x\right)=1\ge 0$. Suppose
that $\mathcal{T}\left(\hat\varphi_{n}(x),x\right)\ge 0$.

Let be $c_n=n(n+1)/4$. Consider the expression
$$
\aligned
&\mathcal{T}\left(\hat\varphi_{n+1}(x),x\right)-c_n\mathcal{T}\left(\hat\varphi_{n}(x),x\right)\\
&= \hat\varphi^2_{n+1}(x)
 - \hat\varphi_{n}(x)\hat\varphi_{n+2}(x) -c_n\left(\hat\varphi^2_{n}(x)-\hat\varphi_{n-1}(x)\hat\varphi_{n+1}(x)\right) \\
& =\hat\varphi^2_{n+1}(x)
 - \hat\varphi_{n}(x)\left(x\hat\varphi_{n+1}(x) - c_{n+1}\hat\varphi_{n}(x)\right)-c_n\hat\varphi_{n-1}(x)\hat\varphi_{n+1}(x) - c_n
 \hat\varphi^2_{n}(x).
 \endaligned
$$
Applying the recurrence relation (\ref{ModMLPRecMonic}), the last
expression reduces to
$$
\mathcal{T}\left(\hat\varphi_{n+1}(x),x\right)-c_n\mathcal{T}\left(\hat\varphi_{n}(x),x\right)
= \left(c_{n+1} - c_n\right)\hat\varphi^2_{n}(x) = \frac{n+1}2
\hat\varphi^2_{n}(x)\ge 0.
$$
We conclude that
$\mathcal{T}\left(\hat\varphi_{n+1}(x),x\right)\ge
c_n\mathcal{T}\left(\hat\varphi_{n}(x),x\right)\ge 0$. $\Box$

\section{Fourier transform}

Let us remind that the Fourier transform is defined by
\begin{equation}\label{fouriertransform}
\mathfrak{F}\left[f(t)\right] = \frac1{\sqrt{2\pi}}\int_{\mathbb
R} f(t)\;e^{{\rm i}st}\;{\rm d}t\;=F(s).
\end{equation}
The Fourier transform of the first members of the sequence
$\{\hat\varphi_n(t)w(t)\}_{n\in\mathbb N_0}$, where $w(t) =
\frac{t}{\sinh (\pi t)}$,  is:
$$
\mathfrak{F}\left[\hat\varphi_{0}w\right]  = \frac1{\sqrt{2\pi}}\; \frac{1}{1+\cosh
s}=\frac1{2}\sqrt{\frac2{\pi}}\; \frac{\sinh^2 (s/2)}{\sinh^2 s},\
\mathfrak{F}\left[\hat\varphi_{1}w\right] = 2\;{\rm i}\;\sqrt{\frac{2}{\pi}}\; \frac{\sinh^4
(s/2)}{\sinh^3 s}.
$$
\begin{thm}
The Fourier transform of the sequence
$\{\hat\varphi_n(t)w(t)\}_{n\in\mathbb N_0}$, where $w(t) =
\frac{t}{\sinh (\pi t)}$,  is:
\begin{equation}\label{Furije}
\Phi_{n}(s) = \mathfrak{F}\left[\hat\varphi_{n}w\right]  = {\rm
i}^n(n+1)!\sqrt{\frac{2}{\pi}}\; \frac{\sinh^{2n+2}
(s/2)}{\sinh^{n+2} s}\qquad (n\in\mathbb N_0).
\end{equation}
\end{thm}
{\it Proof}. We will apply the mathematical induction. It is
obviously true for $n=1$ and a few others. Suppose that it is
valid for every $k\le n$. We will multiply the recurrence relation
(\ref{ModMLPRecMonic}) with $w(t)$ and apply the Fourier transform
on it:
\begin{equation}\label{Furije2}
\mathfrak{F}\left[\hat\varphi_{n+1}(t)w(t)\right] =
\mathfrak{F}\left[t\hat\varphi_{n}(t)w(t)\right] -
\frac{n(n+1)}{4}\
\mathfrak{F}\left[\hat\varphi_{n-1}(t)w(t)\right].
\end{equation}
The Fourier transform (\ref{fouriertransform}) has the property
$$
\mathfrak{F}\left[t^m g(t)\right](s) = (-{\rm i})^m
\frac{d^m}{ds^m}\mathfrak{F}\left[g(t)\right](s) \quad
(m\in\mathbb N).
$$
It able us to write
$$
\mathfrak{F}\left[t \varphi_{n}(t)w(t)\right] = (-{\rm
i})\frac{d}{ds} \mathfrak{F}\left[\varphi_{n}(t)w(t)\right](s) =
-{\rm i}\; \Phi_{n}'(s).
$$
Hence the relation (\ref{Furije2}) obtains the form
$$
\Phi_{n+1}(s) = -{\rm i}\; \hat\Phi'_{n}(s) - \frac{n(n+1)}{4}\
\Phi_{n-1}(s).
$$
Deriving (\ref{Furije}), we find
$$
\Phi'_{n}(s) = {\rm i}^n(n+1)!\sqrt{\frac{2}{\pi}}\;
\bigl(n-2\sinh^2(s/2)\bigr) \frac{\sinh^{2n+2} (s/2)}{\sinh^{n+3}
s}\;.
$$
Finally, using the assumed expression for $\Phi_{n-1}(s)$, we see
that $\Phi_{n+1}(s)$ satisfies the relation (\ref{Furije}),
wherefrom the statement follows. $\Box$

\begin{remark}\rm
Having in mind relations between hyperbolic functions, relation
(\ref{Furije}) can be rewritten in the form
$$
\Phi_{n}(s) = \mathfrak{F}\left[\hat\varphi_{n}w\right]  =
\dfrac{{\rm i}^n(n+1)!}{2^n\;\sqrt{2\pi}}\;
\tanh^n(s/2)\left(1-\tanh^2(s/2)\right)\qquad (n\in\mathbb N_0).
$$
\end{remark}

\section{Differential properties}

The exponential generating function of sequence
$\{\hat\varphi_n(x)\}_{n\in\mathbb N_0}$ is (see \cite{StaMarRaj})
\begin{equation}\label{GenFun222}
\hat{\mathcal G}(t,x)=\frac{4\exp\bigl(2x\;\arctan
(t/2)\bigr)}{t^2+4}=\sum_{n=0}^\infty
\hat\varphi_{n}(x)\frac{t^n}{n!}.
\end{equation}
\begin{thm}
The exponential generating function $\hat{\mathcal G}(t,x)$ has
the property
\begin{equation*}%\label{GenFun223}
\hat{\mathcal G}(t,x)\; \frac{\partial^2}{\partial
x^2}\hat{\mathcal G}(t,x) = \left(\frac{\partial}{\partial
x}\hat{\mathcal G}(t,x)\right)^2.
\end{equation*}
\end{thm}
\noindent{\it Proof.} This follows directly from relations
$$
\frac{\partial}{\partial x}\hat{\mathcal G}(t,x) = 2\arctan
\frac{t}{2}\cdot\hat{\mathcal G}(t,x) \;, \quad
\frac{\partial^2}{\partial x^2}\hat{\mathcal G}(t,x) = 4\arctan^2
\frac{t}{2} \cdot \hat{\mathcal G}(t,x) \;.\;\Box
$$
\begin{cor}
The sequence $\{\hat\varphi_n(x)\}_{n\in\mathbb N_0}$ satisfies
the recurrence-differential equation
\begin{equation*}%\label{GenFun224}
\sum_{k=0}^n \frac{\hat\varphi''_k(x)\hat\varphi_{n-k}(x) -
\hat\varphi'_k(x)\hat\varphi'_{n-k}(x)}{k!(n-k)!} = 0\quad
(n\in\mathbb N).
\end{equation*}
\end{cor}

%H. Youn, Y. Yang, Differential Equation and Recursive Formulas of Sheffer Polynomial Sequences,
% International Scholarly Research Network ISRN Discrete Mathematics Volume 2011, Article ID 476462, 16 pages doi:10.5402/2011/476462
\begin{thm}
Any polynomial $\hat\varphi_n(x)$ satisfies $n^{\rm th}$ order
differential equation of the form
\begin{equation}\label{DE-Mit-Leff-1}
\sum_{k=1}^n
(\alpha_k+\beta_kx)\frac{\hat\varphi^{(k)}_{n}(x)}{k!}
-n\hat\varphi_{n}(x)=0,
\end{equation}
where  $\alpha_k = \cos\frac{k\pi}2$ and $\beta_k =
\sin\frac{k\pi}2$.
\end{thm}
{\it Proof}. The sequence $\{\hat\varphi_n(x)\}_{n\in\mathbb N_0}$
is a Sheffer sequence since its generating function has the form
\begin{equation*}%\label{Sheffer}
\hat{\mathcal G}(t,x) = \frac1{g\left(\hat f(t)\right)}{\rm
e}^{x\hat f(t)},
\end{equation*}
where $\hat f(t)=2\arctan \frac{t}2$. It is the compositional
inverse of $f(t)=2\tan \frac{t}2$. Also, here is $ g\left(\hat
f(t)\right) = 1+\frac{t^2}4. $ Hence $g(t) = \cos^2\frac{t}2$.

According to \cite{Youn}, $\{\hat\varphi_n(x)\}_{n\in\mathbb N_0}$
satisfies the differential equation of the form
(\ref{DE-Mit-Leff-1}), where
$$
\beta_k =
\left(\frac{f(t)}{f'(t)}\right)^{(k)}\Big|_{t=0}\;,\qquad \alpha_k
= \left(-\frac{f(t)}{f'(t)}\cdot
\frac{g'(t)}{g(t)}\right)^{(k)}\Big|_{t=0}\;.
$$
Hence
$$
\beta_k =  \left(\sin t\right)^{(k)}\Big|_{t=0} =
\sin\frac{k\pi}2\;,\qquad
\alpha_k = \left(1-\cos t\right)^{(k)}\Big|_{t=0} =
\cos\frac{k\pi}2\;.\Box
$$
\begin{example}\rm The polynomial
$$
\hat\varphi_{4}(x) = x^4-5x^2+\frac32,
$$
satisfies the following differential equation:
\begin{equation}\label{DE-Mit-Leff-11}
\frac1{24}\hat\varphi^{(4)}_{4}(x)
-\frac1{6}x\hat\varphi^{(3)}_{4}(x) -\frac1{2}\hat\varphi''_{4}(x)
+x\hat\varphi'_{4}(x)-4\hat\varphi_{4}(x)= 0.
\end{equation}
\end{example}
%Thm 2.7
\begin{thm}
Any polynomial $\hat\varphi_n(x)$ satisfies differential equation
of the form
\begin{equation*}%\label{DE-Mit-Leff-2}
\left( \cos D + x\sin D - (n+1)I\right) \hat\varphi_n(x) = 0
\qquad \left(D=\frac{{\rm d}}{{\rm d}x}\right).
\end{equation*}
The polynomial $\hat\varphi_n(x)$ is the eigenfunction of the
operator $\mathcal{F} = \cos D + x\sin D - I$ with the eigenvalue
$n$.
\end{thm}
{\it Proof}. Let $I$ be the identity operator. Since
$$
\alpha_k =  \cos\frac{k\pi}2 = \frac{{\rm i}^k+(-{\rm
i})^k}2\;,\qquad \beta_k = \sin\frac{k\pi}2 = \frac{{\rm
i}^k-(-{\rm i})^k}{2{\rm i}}\;,
$$
we can write (\ref{DE-Mit-Leff-1}) in the form
\begin{equation}\label{DE-Mit-Leff-3}
\left(\sum_{k=1}^n \left(\frac{{\rm i}^k+(-{\rm i})^k}2 +
x\frac{{\rm i}^k-(-{\rm i})^k}{2{\rm i}}\right)\frac{D^k}{k!}
-nI\right)\hat\varphi_{n}(x)=0.
\end{equation}
Since $D^m \hat\varphi_{n}(x)\equiv 0$ for every $m>n$, we can
write
$$
\sum_{k=1}^n {\rm i}^k\frac{D^k}{k!}\hat\varphi_{n}(x) =
\sum_{k=1}^\infty \frac{({\rm i}D)^k}{k!}\hat\varphi_{n}(x) =
\left({\rm e}^{{\rm i}D}-I\right)\hat\varphi_{n}(x).
$$
Hence the formula (\ref{DE-Mit-Leff-3}) becomes
$$
\left(\frac12 \left( {\rm e}^{{\rm i}D}+ {\rm e}^{-{\rm
i}D}-2I\right) +\frac{x}{2{\rm i}}\left( {\rm e}^{{\rm i}D} - {\rm
e}^{-{\rm i}D}\right) - nI\right)\hat\varphi_{n}(x)=0.
$$
The statement follows from the Euler identity for the complex
functions. $\Box$
\begin{example}\rm
Since
$$
\cos D = \sum_{k=0}^\infty (-1)^k\frac{D^{2k}}{(2k)!} \;,\qquad
\sin D = \sum_{k=0}^\infty (-1)^k\frac{D^{2k+1}}{(2k+1)!}\;,
$$
the polynomial $\hat\varphi_{4}(x)$ satisfies
$$
\left( \left(I-\frac{D^{2}}{2} + \frac{D^{4}}{4!} \right)+ x\left(
D - \frac{D^{3}}{3!}\right)-4I\right)\hat\varphi_{4}(x)= 0,
$$
what is the same as (\ref{DE-Mit-Leff-11}).
\end{example}
\begin{thm}
The sequences $\{\hat \varphi_n(x)\}_{n\in\mathbb N_0}$ and
$\{\varphi_n(x)\}_{n\in\mathbb N_0}$ have the fo\-llow\-ing
differential properties:
\begin{equation}\label{DiffPro_NN}
\hat\varphi'_{n+1}(x) = \sum_{k=0}^{[n/2]}
(-1)^{k}\binom{n+1}{2k+1}\frac{(2k)!}{2^{2k}}\hat\varphi_{n-2k}(x),
\end{equation}
\begin{equation}\label{DiffPro_N}
\varphi'_n(x) = 2\sum_{k=0}^{[n/2]} \frac{(-1)^{k}}{2k+1}
\varphi_{n-2k}(x).
\end{equation}
\end{thm}

\noindent{\it Proof.} By differentiation the generating function
(\ref{GenFun222}) over $x$, we get
\begin{equation*}
\sum_{n=0}^\infty \hat\varphi'_{n}(x)\frac{t^n}{n!}=
\frac{4\exp\bigl(2x\ \arctan (t/2)\bigr)}{t^2+4}\ 2\arctan (t/2).
\end{equation*}
Knowing that $\hat\varphi'_0(t)=0$ and using the expansion
$$
2\arctan\dfrac{t}{2}=\sum_{k=0}^\infty
\dfrac{(-1)^k}{4^k(2k+1)}t^{2k+1}\qquad
\left(\Big|\dfrac{t^2}{2}\Big|<1\right),
$$
we have
$$
\sum_{n=1}^\infty \hat\varphi'_{n}(x)\frac{t^n}{n!}=
\left(\sum_{n=0}^\infty
\hat\varphi_{n}(x)\frac{t^n}{n!}\right)\left(\sum_{k=0}^\infty
\dfrac{(-1)^k}{4^k(2k+1)}t^{2k+1}\right).
$$
Hence
$$
\sum_{n=0}^\infty \frac{\hat\varphi'_{n+1}(x)}{(n+1)!}t^{n+1}=
\sum_{n=0}^\infty\sum_{k=0}^{\infty}\dfrac{\hat\varphi_{n}(x)}{n!}
\dfrac{(-1)^k}{4^k(2k+1)}t^{n+2k+1},
$$
i.e.,
$$
t\sum_{n=0}^\infty \frac{\hat\varphi'_{n+1}(x)}{(n+1)!}t^{n}=
t\sum_{n=0}^\infty\sum_{k=0}^{[n/2]}\dfrac{\hat\varphi_{n-2k}(x)}{(n-2k)!}
\dfrac{(-1)^k}{4^k(2k+1)}t^{n}.
$$
Comparing the coefficients by $t^n\ (n\in\mathbb N)$, we find
\begin{equation}\label{pom}
\frac{\hat\varphi'_{n+1}(x)}{(n+1)!}
=\sum_{k=0}^{[n/2]}\dfrac{\hat\varphi_{n-2k}(x)}{(n-2k)!}
\dfrac{(-1)^k}{4^k(2k+1)}.
\end{equation}
By rearrangement of summands, we have formula (\ref{DiffPro_NN}).

Formula (\ref{DiffPro_N}) can be obtained by (\ref{pom}) and
(\ref{ModMLPMonic}).$\quad\Box$

\section{Quasi-monomiality}

%Youssef Ben Cheikh, Some results on quasi-monomiality, Applied Mathematics and Computation 141 (2003) 63–76

% H. Chaggara, W. Koepf, Duplication coefficients via generating
% functions, Complex Variables and Elliptic Equations, Vol. 52, No. 6, (2007) 537–549
% Dodati Y. Ben Cheikh a,., H. Chaggarab, Connection coefficients between Boas–Buck
%polynomial sets J. Math. Anal. Appl. 319 (2006) 665–689 *Cheikn.pdf

According to  \cite{Cheikh} and \cite{Chaggara}, the exponential
generating function ${\mathcal G}(t,x)$ is of the Boas-Buck type
if
\begin{equation*}%\label{Boas-Buck}
{\mathcal G}(t,x) = A(t)B(xC(t)),
\end{equation*}
where
$$
B^{(k)}(0)\ne  0\quad (\forall k\in\mathbb N), \quad A(0)C'(0)\ne
0, \quad C(0)= 0.
$$
Considering the exponential generating function of the reduced
Mittag-Leffler polynomials $\{\hat \varphi_n(x)\}$, we can denote
with
$$
A_{\hat\varphi}(t)=\frac{4}{4+t^2},\quad B_{\hat\varphi}(t)={\rm
e}^t, \quad C_{\hat\varphi}(t)=2\arctan \frac{t}2.
$$
Here, it is
$$
C^{-1}_{\hat\varphi}(t) = 2\tan\frac{t}{2}\;.
$$
\begin{thm}
The sequence $\hat\varphi_n(x)$ is quasi-monomial under the
lowering operator $ \mathcal{L}_x = 2\tan\left(D_x/2\right) $,
i.e.
\begin{equation*}%\label{quasi-mon-1}
\mathcal{L}_x \hat\varphi_{n}(x) = n\;\hat\varphi_{n-1}(x)\qquad
(n\in\mathbb N).
\end{equation*}
\end{thm}
\noindent{\it Proof.} We start with the Taylor series
$$
\tan x = \sum_{k=1}^\infty \theta_{2k-1} x^{2k-1}, \quad  {\rm
where}\quad \theta_{2k-1} =
(-1)^{k-1}4^k(4^k-1)\frac{B_{2k}}{(2k)!}\;.
$$
Here, $B_n$ is the $n$-the Bernoulli number. Since
$$
D_x^{2k-1} \hat{\mathcal G}(t,x) =
\frac{4}{t^2+4}\exp\Bigl(2x\,\arctan \frac{t}2\Bigr)\;
\bigl(2\,\arctan \frac{t}2\bigr)^{2k-1},
$$
we have
$$
\aligned \mathcal{L}_x \hat{\mathcal G}(t,x) &= 2\sum_{k=1}^\infty
\theta_{2k-1} \left(\frac{D_x}2\right)^{2k-1} \hat{\mathcal
G}(t,x) \\
&= \frac{8}{t^2+4}\exp\Bigl(2x\,\arctan \frac{t}2\Bigr)
\sum_{k=1}^\infty \theta_{2k-1} \Bigl(\arctan
\frac{t}2\Bigr)^{2k-1},
\endaligned
$$
wherefrom
$$
\mathcal{L}_x \hat{\mathcal G}(t,x) = t\, \hat{\mathcal G}(t,x).
$$
Since
$$
\mathcal{L}_x \hat{\mathcal G}(t,x) = \mathcal{L}_x \left(
 \sum_{n=0}^\infty \hat\varphi_{n}(x)\frac{t^n}{n!}\right)
=\sum_{n=0}^\infty \mathcal{L}_x
\hat\varphi_{n}(x)\frac{t^n}{n!}\;,
$$
and
$$
t\; \hat{\mathcal G}(t,x) = t\left(
 \sum_{n=0}^\infty \hat\varphi_{n}(x)\frac{t^n}{n!}\right)
= \sum_{n=1}^\infty n\;\hat\varphi_{n-1}(x)\frac{t^{n}}{n!},
$$
we have the statement proven. $\Box$

%\begin{\bf Acknowledgement}
{\bf Acknowledgement}. This research was financially supported by
the Ministry of Education, Science and Technological Development
of the Republic of Serbia.
%\end{acknowledgement}

\noindent {\sc Predrag M. Rajkovi\'c}\\
Faculty of Mechanical Engineering, University of Ni\v s\\
A. Medvedeva 14,18 000 Ni\v s, Serbia \\
{\it E-mail}:\ predrag.rajkovic@masfak.ni.ac.rs \\
[8mm]\\
{\sc Sladjana D. Marinkovi\'c}\\
Faculty of Electronic Engineering, University of Ni\v s\\
A. Medvedeva 12,18 000 Ni\v s, Serbia \\
{\it E-mail}:\ sladjana.marinkovic@elfak.ni.ac.rs \\
[2mm]\\
{\sc Miomir S. Stankovi\'c}\\
The Mathematical Institute of  SASA, Belgrade, Serbia, \\
Kneza Mihaila 36, 11 000 Belgrade, Serbia\\
{\it E-mail}:\ miomir.stankovic@gmail.com \\
[2mm]\\
{\sc Marko D. Petkovi\'c}\\
Faculty of Mathematics and Sciences, University of Ni\v s\\
Vi\v segradska bb,18 000 Ni\v s, Serbia \\
{\it E-mail}:\ dexterofnis@gmail.com

\end{document}